\documentclass[12pt,oneside]{amsart}
\usepackage{amssymb, amsmath, amsthm}

\usepackage{epsfig}
\usepackage{epstopdf}
\usepackage{graphicx}
\usepackage{color}
\usepackage{enumerate}
\usepackage{amsrefs}

\usepackage[normalem]{ulem}

\usepackage[normalem]{ulem}

\usepackage{pinlabel}

\theoremstyle{plain}
\newtheorem{theorem}{Theorem}[section]

\newtheorem*{theorem*}{Theorem}
\newtheorem*{theorem-DisjtSS}{Theorem \ref{Thm: Disjt SS}}
\newtheorem*{theorem-EssentialTorus}{Theorem \ref{Thm: Essential Torus}}
\newtheorem*{cor-ScharlWu}{Corollary \ref{Cor-ScharlWu}}
\newtheorem*{corollary-MSC 2}{Corollary \ref{Cor: MSC 2}}
\newtheorem*{theorem-Main A}{Theorem \ref{Thm: Main A}}
\newtheorem*{theorem-Main B}{Theorem \ref{Thm: Main Thm B}}
\newtheorem*{Cor-Unknotting}{Theorem \ref{Cor: Prime Unknotting 1}}
\newtheorem*{Cor-genusbandsum}{Theorem \ref{Thm: Genus superadd}}
\newtheorem*{Cor-bandsumscc}{Corollary \ref{Cor: Band Sums CC}}

\newtheorem{proposition}[theorem]{Proposition}

\newtheorem{lemma}[theorem]{Lemma}
\newtheorem{definition}[theorem]{Definition}

\theoremstyle{definition}

\newtheorem*{Claim}{Claim}

      \makeatletter
      \def\@setcopyright{}
      \def\serieslogo@{}
      \makeatother

\begin{document}

   \title[Self-replicating 3-manifolds]{Self-replicating 3-manifolds}
   \author{Ryan Blair and Ricky Lee}
   \email{}
   \thanks{}

\begin{abstract}
In this paper we explore the topological properties of self-replicating, 3-dimensional manifolds, which are modeled by  idempotents in the (2+1)-cobordism category. We give a classification theorem for all such idempotents. Additionally, we characterize biologically interesting ways in which self-replicating 3-manifolds can embed in $\mathbb{R}^3$.
\end{abstract}

\maketitle
\date{\today}

\section{Introduction}
The motivation of this paper is the study of self-replicating 3-manifolds. Intuitively, a 3-manifold $M$ is self-replicating if it contains a surface $F$ such that cutting along $F$ yields two components, each homeomorphic to $M$. Kauffman previously suggested that the natural model for general self-replication is that of idempotents in appropriately topological categories  \cite{Kauff04}. Previous classification of idempotents in topological categories include the Temperley-Lieb category \cite{Abramsky07} and the tangle category \cite{BS19}.  In this paper we model self-replicating 3-manifolds as idempotents in the $(2+1)$-cobordism category and classify all such morphisms. We also explore questions of the embedability of self-replicating 3-manifolds in $\mathbb{R}^3$.

An \emph{idempotent} of a category is a morphism that is idempotent with respect to composition, i.e.\ a morphism $M$ such that $M = M\circ M$. Idempotents have applications to the theory of quantum observation \cite{Selinger08}, self-replication in biology 
\cite{Kauff04}, and algebraic structures associated to quantum theory \cite{HaqKauffman}. An idempotent $M$ \emph{splits} if there are morphisms $P$ and $Q$ such that $M=P\circ Q$ and $Q\circ P$ is an identity morphism. Any morphism $M$ that splits is idempotent (if $Q\circ P$ is an identity, then $(P\circ Q)\circ (P\circ Q) = P\circ (Q\circ P)\circ Q = P\circ Q$).  However, in many categories, not all idempotents split. Let $\mathcal{C}$ denote the $(2+1)$-cobordism category. In this paper, the objects in $\mathcal{C}$ are compact, orientable surfaces. The set of morphisms of $\mathcal{C}$ from $F_1$ to $F_2$ is denoted $Mor(F_1,F_2)$ and consists of compact, orientable 3-manifolds $M$ for which $F_1 \cup F_2$ naturally embeds in $\partial M$.

\begin{theorem}\label{main1}
If $M\in Mor(G,G)$ is an idempotent such that such that $M$ and $G$ are connected as manifolds, then $M$ splits.

\end{theorem}

The above theorem establishes the first goal of this paper, a classification of self-replicating 3-manifolds. Idempotents that split are particularly simple since they are in one-to-one correspondence with decompositions of the identity morphisms. It was previously shown that all idempotents in the category of unoriented tangles up to isotopy split \cite{BS19}. Our proof of the above theorem is inspired by the strategy employed in that paper. However, there are several key technical differences.

The second goal of this paper is to explore the concept of realizability of self-replicating 3-manifolds in $\mathbb{R}^3$. As we will show in Proposition \ref{AllEmbed}, all idempotents embed in $\mathbb{R}^3$. So, we instead ask about the more restrictive question of effective embedding. Given a decomposition of a morphism $M$ of $\mathcal{C}$ as $M=M_1 \circ M_2$, the \emph{decomposing surface} $F$ is the properly embedded surface in $M$ corresponding to the codomain of $M_2$ and the domain of $M_1$. An idempotent 3-manifold $M\in Mor(F,F)$ has an \emph{effective} embedding into $\mathbb{R}^3$ if the image of $M$ in $\mathbb{R}^3$ can be surgered along the decomposing surface $F$ corresponding to $M=M\circ M$ to produce two embeddings of $M$, denoted $M_1$ and $M_2$, such that there is an embedded $2$-sphere in $\mathbb{R}^3$ separating $M_1$ from $M_2$. Note that an effective embedding of an idempotent is meant to model the possibility of biological dispersal, movement that has the potential to lead to gene flow. In particular, the embeddings $M_1$ and $M_2$ in the definition of effective embedding are unlinked in $\mathbb{R}^3$. 

\begin{figure}[h!]
\begin{picture}(255,410)
\put(1,1){\includegraphics[scale=.5]{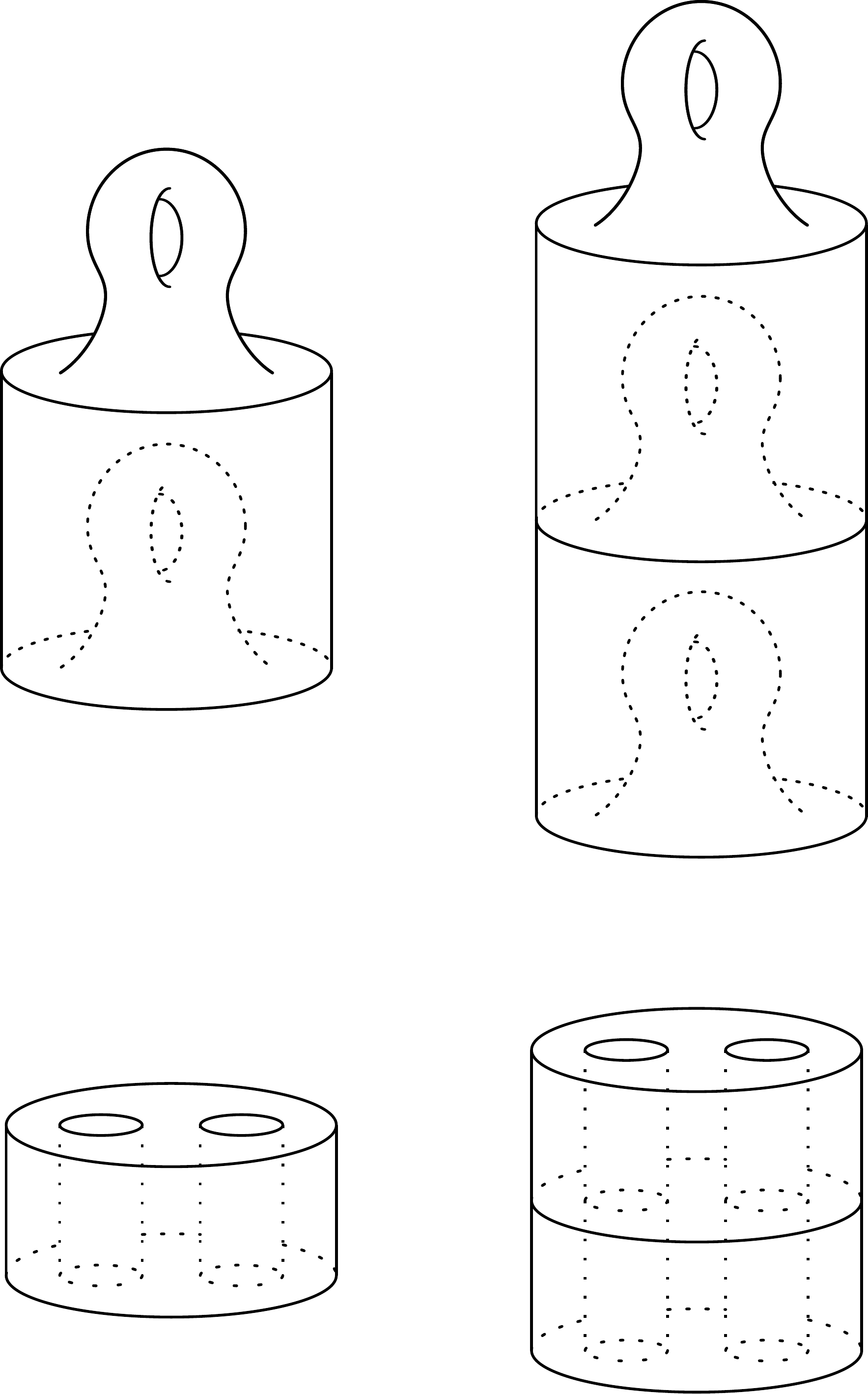}}
\put(123,55){$\cong$}
\put(123,270){$\cong$}
\end{picture}
\caption{Two representations of the genus 2 handlebody as an idempotent in the $(2+1)$-cobordism category. The top representation is not an effective embedding. The bottom representation is a trivial idempotent and an effective embedding.}
\label{G2example}
\end{figure}

We say an idempotent in $\mathcal{C}$ is \emph{trivial} if it is the trivial morphism on some compact surface $F$. In particular, a trivial morphism is homeomorphic to $F\times I$. We call a surface \emph{planar}, if it embeds in $S^2$. All trivial morphisms on planar surfaces have effective embeddings into $\mathbb{R}^3$. However, the trivial morphism on the torus does not. See Figure \ref{G2example} for an example of a trivial idempotent with an effective embedding and an example of a non-trivial idempotent realized by an embedding that is not effective.

The following theorem gives a characterization of self-replicating 3-manifolds with effective embeddings into $\mathbb{R}^3$.

\begin{theorem}\label{main2}
Suppose $G$ is a connected compact orientable surface and $M$ is a connected idempotent with $M\in Mor(G,G)$. The 3-manifold,  $M$, effectively embeds into $\mathbb{R}^3$ if and only if $M$ is a trivial morphism and $G$ is planar. 
\end{theorem}

The paper is organized as follows. In Section \ref{Sec:Prelim} we give a rigorous definition of the $(2+1)$-cobordism category and establish the 3-manifold machinery needed for the subsequent proofs. In Section \ref{Sec:Split} we prove Theorem \ref{main1}. In Section \ref{Sec:Effective}, we prove Theorem \ref{main2}.

\section{Preliminaries}\label{Sec:Prelim}

\subsection{The $(2+1)$ Cobordism Category}

We begin this section with a description of the $(2+1)$-cobordism category. For additional details on this category see \cite{CR18} and \cite{Milnor1}. A \emph{triad} is a triple $(M,F,G)$ where $M$ is a smooth compact 3-manifold and $F$ and $G$ are smooth, compact, orientable 2-manifolds such that $\partial M = (F \amalg G) \cup X$ with $X \cong (\partial F) \times I$.  When $\partial F$ is empty, $\partial M = F \amalg G$. Otherwise, $M$ is
a smooth manifold with corners along $\partial F$ and $\partial G$. However, we will often suppress the manifold with corners structure on $M$ since we will often be content to classify $M$ up to homeomorphism.
Given fixed surfaces $F$ and $G$, a cobordism from $F$ to $G$ is a triple $(M, j_F, j_G)$, where $M$ is a compact
smooth 3-manifold and $j_F : F \rightarrow \partial M$, $j_G : G \rightarrow \partial M$ are embeddings such
that $(M, j_{F}(F), j_{G}(G))$ is a triad. When appropriate, we will sometimes suppress the additional structure and refer to the cobordism $(M, j_F, j_G)$ as $M$.

Two cobordisms $(M_1, j_{F_1}, j_{G_1})$ and $(M_2, j_{F_2}, j_{G_2})$ are equivalent if there is a diffeomorphism $h:M_1\rightarrow M_2$ such that $h \circ j_{F_1}=j_{F_2}$ and $h \circ j_{G_1}=j_{G_2}$. The set of equivalence classes of cobordisms from $F$ to $G$ is denoted by $Mor(F,G)$. In this more precise definition the trivial cobordism form $F$ to $G$ is $(F\times I, j^0_{F},j^1_{F})$ where $j^i_F:F\rightarrow F\times I$ is the standard inclusion map such that $j^i_F(F)=F\times \{i\}$. 

Let $M\in Mor(F,G)$ and $M'\in Mor(G,H)$ be represented by cobordisms $(M, j_F,j_G)$ and $(M', j'_G,j'_H)$. The topological manifold $M\cup_{j'_G\circ (j_G)^{-1}}M'$ admits a smooth structure compatible with those of $M$ and $M'$, giving rise to a smooth manifold $M\circ M'$, and $(M\circ M',j_F, j'_H)$ represents a well-defined class $M\circ M' \in Mor(F,H)$. Then, an idempotent is an equivalence class of cobordism $M \in Mor(F,F)$ such that $M=M\circ M$.

\subsection{Decomposing Surfaces}

	\begin{defn}
		Let $M\in Mor(F,G)$ be represented by the cobordism $(M, j_F,j_G)$. A surface $H\subset M$ is a \emph{decomposing surface} for $M$ if:
		\begin{enumerate}
			\item $H$ is a smooth, compact, orientable 2-manifold.
			\item $j_H:H\rightarrow M$ is a proper smooth embedding.
			\item If $\partial F \neq \emptyset$, then for each annular component $A_i$ of $\partial M\setminus (j_{F}(F) \amalg j_{G}(G))$ there exists a unique connected component $\alpha_i$ in $\partial H$ such that $\alpha_i$ is an essential loop in $A_i$. 
			\item The exterior of $H$ in $M$ is  $A \coprod B$, where $A$ and $B$ are  3-manifolds such that $F \subset \partial A$ and $G \subset \partial B$. Moreover, $A\in Mor(F,H)$, $B\in Mor(H,G)$, and $B\circ A$ is equivalent to  $M$ as an element of $Mor(F,G)$.
		\end{enumerate}
	\end{defn}
	
	Note that a consequence of this definition is that for any decomposing surface $H$ for $M\in Mor(F,G)$, each of $\partial H$, $\partial F$ and $\partial G$ contain the same number of components.

	\begin{defn}
		A decomposing surface $H$ for $M\in Mor(F,G)$ is \emph{minimal} if there is no decomposing surface $H'$ for $M$ with $-\chi(H')< -\chi(H)$.
	\end{defn}

\subsection{Essential Surfaces}

A surface $F$ properly embedded in a 3-manifold $M$ is \emph{boundary-parallel} if there is an isotopy of $F$ in $M$ which fixes $\partial F$ and takes $F$ to a subsurface contained in $\partial M$. Otherwise, we say $F$ is \emph{non-boundary parallel}. A loop $\gamma$ embedded in a surface $F$ is \emph{essential} if it does not bound a disk in $F$. A surface $F$ is \emph{compressible} in $M$ if $F$ is a $2$-sphere bounding a $3$-ball or if
there exists a disk $D$ embedded in $M$ such that $D \cap F=\partial D$ and $\partial D$ is essential in $F$. Such a disk is called a \emph{compressing disk}.  Otherwise, we say $F$ is \emph{incompressible}. A surface in $M$ is \emph{essential} if it is incompressible and non-boundary parallel.

Given a properly embedded surface $F$ in a 3-manifold $M$, we can \emph{compress} $F$ along a compressing disk $D$ to form a new properly embedded surface $F^*$.  Let $D^2\times I$ be a small fibered neighborhood of $D$ in $M$ such that $D=D^2\times\{\frac{1}{2}\}$ and $\partial(D^2) \times I$ is an embedded annulus in $F$. Then we define $F^*$ to be the surface isotopic to $(F\setminus (\partial(D^2) \times I)) \cup (D^2\times \{0,1\})$.

The number of non-isotopic, disjoint, essential surfaces properly embedded in a compact 3-manifold is bounded due to the following classical result.
	
	\begin{theorem}\label{Parallel} [Page 49 of \cite{Jaco80}]
		Let $M$ be a compact 3-manifold. If $\{F_1, \ldots ,F_n\}$ is a collection of pairwise disjoint, incompressible surfaces in $M$ so that for some integer $\chi_0$, $\chi(F_i)>\chi_0$, $1\leq i\leq n$, then there is an integer $N_0(M, \chi_0)$ such that either $n< N_0(M, \chi_0)$, some $F_i$ is an annulus, or a disk parallel into $\partial M$, or for some $i \neq j$, $F_i$ is parallel to $F_j$ in $M$.
	\end{theorem}
	
	Note that the conclusion of $F_i$ being parallel to $F_j$ in the above theorem implies that $F_i\cong F_j \cong F$ and that $F_i \cup F_j$ bounds a 3-manifold homeomorphic to $F\times I$ in $M$ such that $F_i=F\times\{0\}$ and $F_j=F\times\{1\}$.
	
\section{Idempotents Split}\label{Sec:Split}
	
	\begin{lemma}\label{Essential}
		Let $M\in Mor(G,G)$ be a nontrivial idempotent such that $M$ and $G$ are connected as manifolds. If $F$ is a minimal decomposing surface for $M$, then $F$ is essential.
	\end{lemma}

	\begin{proof}
		Suppose $F$ is compressible. Then we can compress $F$ once to obtain a surface $F^*$. The surface $F^*$ may have more connected components than $F$. However, $\partial F = \partial F^*$ and  $\chi(F^*)=\chi(F)+2 > \chi(F)$. Moreover, $F^*$ continues to be a separating surface in $M$ with $\partial_+ M$ to one side and $\partial_- M$ to the other. Hence, $F^*$ is a decomposing surface for $M$ with $-\chi(F^*)< -\chi(F)$, contradicting the minimality of $F$. Thus, $F$ is incompressible.

		Now we show $F$ cannot be boundary parallel. We begin by proving the following claim:
	
	\begin{Claim}
		Suppose $M\in Mor(G,G)$ is a nontrivial idempotent. Then either $\partial_+M$ or $\partial_- M$ is compressible in $M$.
	\end{Claim}

	\begin{proof}
		Suppose for contradiction both $\partial_+ M$ and $\partial_- M$ are incompressible. By Theorem \ref{Parallel}, there exists integers $N_0$ and $\chi_0$ such that if $F_1, \ldots ,F_k$ is a collection of pairwise disjoint incompressible surfaces with $k>N_0$ and $\chi(F_i)> \chi_0$, for all $1\leq i\leq k$, then some $F_i$ is a boundary parallel annulus, boundary parallel disk, or there exists $i\neq j$ such that $F_i$ is parallel to $F_j$. Fix such integers $N_0$ and $\chi_0$ with $\chi_0 < \chi(G)$.\par 
		
Since $M$ is an idempotent, $M \cong M\circ M\circ M \ldots M$ where we have composed $M$ with itself $N_0+2$ times. Therefore, we can find $N_0+1$ disjoint decomposing surfaces, $G_1, \ldots ,G_{N_0+1}$ in $M$, each of which is homeomorphic to  $G$. See Figure \ref{schematic} for the case when $N_0=1$. Since each $G_i$ is homeomorphic to $G$, then $\chi_0 < \chi(G_i)$ for all $1\leq i\leq N_0+1$. Since both $\partial_+ M$ and $\partial_- M$ are incompressible, then each $G_i$ is incompressible.

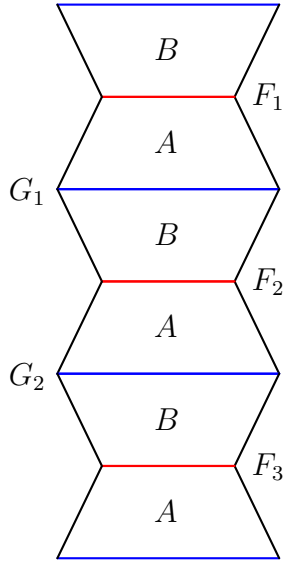
\begin{figure}[h!]
\begin{picture}(85,210)
\put(1,1){\includegraphics[scale=.5]{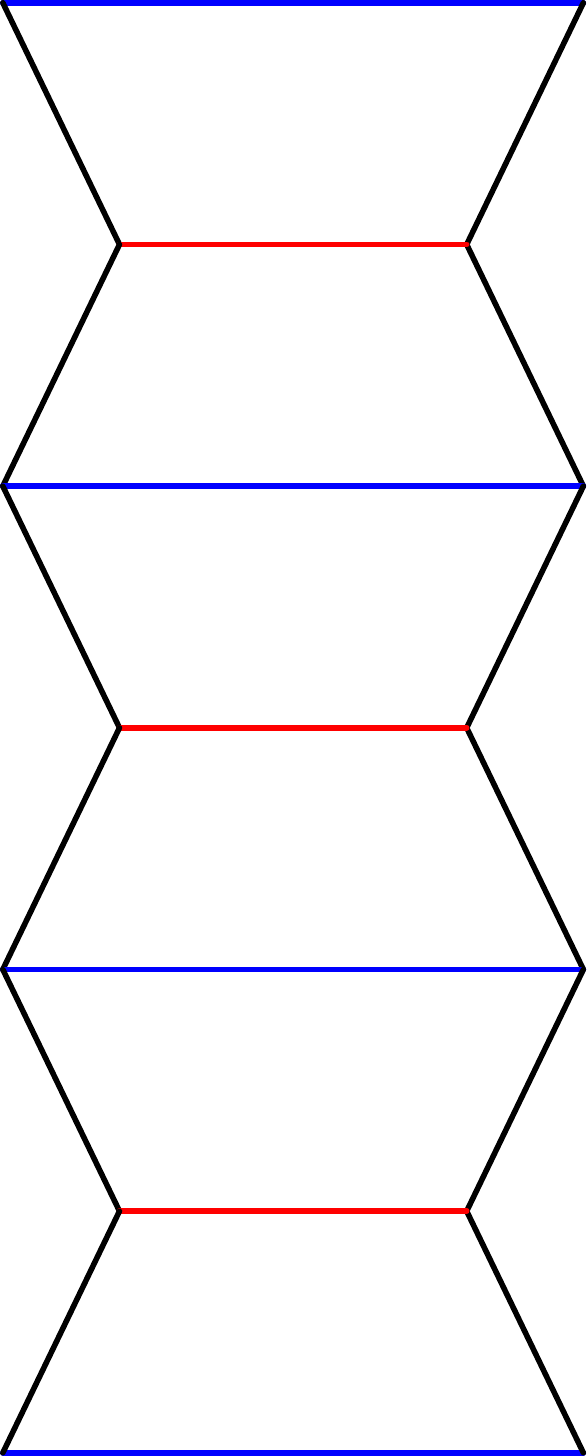}}
\put(38,15){$A$}
\put(38,50){$B$}
\put(38,85){$A$}
\put(38,120){$B$}
\put(38,155){$A$}
\put(38,190){$B$}
\put(-17,67){$G_2$}
\put(-17,137){$G_1$}
\put(75,33){$F_3$}
\put(75,103){$F_2$}
\put(75,173){$F_1$}
\end{picture}
\caption{A schematic decomposition of an idempotent $M$ as $M\circ M \circ M$ where $M=B\circ A$.}
\label{schematic}
\end{figure}
		
		Since the collection of surfaces, $G_1, \ldots ,G_{N_0+1}$ meets the hypothesis for Theorem \ref{Parallel}, then one of the following hold: some $G_i$ is a boundary parallel annulus, some $G_i$ is a boundary parallel disk, or there exists $i\neq j$ such that $G_i$ is parallel to $G_j$. In each of these cases $M$ is of the form $surface\times I$. But $M$ was assumed to by a nontrivial idempotent, so we get our desired contradiction.
	\end{proof}
	
		Suppose $F$ is boundary parallel. Then $F\cong G$, and $\chi(F) = \chi(G)$. By our claim, we can assume one of $\partial_- M$ and $\partial_+ M$ is compressible. Without loss of generality, assume $\partial_+ M$ is compressible.  Let $G\times I$ denote a collar neighborhood of $\partial_+ M$ with $\partial_+ M = G \times \{1\}$. Set $N:= G\times \{0\}$. Then $N$ is a decomposing surface such that $\chi(N)=\chi(\partial_+M)=\chi(F)$. Since $F$ is a minimal decomposing surface, then $N$ is also a minimal decomposing surface. Since $\partial_+M$ is compressible, then $N$ is also compressible. But minimal decomposing surfaces were shown to be incompressible. This gives our desired contradiction. Hence, $F$ is essential.
	\end{proof}

The following is the proof of Theorem \ref{main1}.

\begin{proof}
	Let $F$ be a minimal decomposing surface for $M$. By Lemma \ref{Essential}, $F$ is essential. The surface $F$ decomposes $M$ into cobordisms $A\in Mor(G,F)$, $B\in Mor(F,G)$, such that $M \cong B \circ A$.\par

	By Theorem \ref{Parallel}, there exists integers $N_0$ and $\chi_0$ such that if $F_1, \ldots , F_k$ is a collection of pairwise disjoint incompressible surfaces with $N_0<k$ and $\chi_0 < \chi(F_i)$ for $1\leq i\leq k$, then some $F_i$ is a boundary parallel annulus, boundary parallel disk, or there exists $i\neq j$ such that $F_i$ is parallel to $F_j$. Fix such integers $N_0$ and $\chi_0$.\par 
	
	Since $M$ is an idempotent, $M \cong M\circ M\circ M \ldots M$ where we have composed $M$ with itself $N_0+1$ times. Since $F$ is an essential minimal decomposing surface, we can find $N_0+1$ disjoint essential minimal decomposing surfaces $F_1, \ldots ,F_{N_0+1}$ for $M$, each representing the copy of $F$ in each copy of $M$. The collection of surfaces $F_1, \ldots ,F_{N_0+1}$ decompose $M$ into one copy of $A$, one copy of $B$, and $N_0$ copies of $A\circ B$. See Figure \ref{schematic} for the case when $N_0=2$. \par

	Note that each surface $F_i$ in the collection $F_1, \ldots ,F_{N_0+1}$ is a minimal decomposing surface for $M$ and essential by Lemma \ref{Essential}. In particular, each surface $F_i$ is not a boundary parallel annulus and not a boundary parallel disk. By Theorem \ref{Parallel}, there exist two surfaces $F_i$ and $F_j$ in our collection of essential minimal decomposing surfaces that are parallel. This shows that the cobordism bounded by $F_i$ and $F_j$ in $M$, which is equivalent to $(A\circ B)^l$ for some $1\leq l$, must be equivalent to the trivial cobordism $F\times I$. Thus $A\circ B$ is the trivial morphism in $Mor(F,F)$.
\end{proof}

\section{Effective Embeddings}\label{Sec:Effective}

\begin{proposition}\label{AllEmbed}
	If $M\in Mor(G,G)$ is an idempotent in $\mathcal{C}$ such that such that $M$ and $G$ are connected as manifolds, then $M$  embeds in $\mathbb{R}^3$
\end{proposition}

\begin{proof}
	Let $M\in Mor(G,G)$ be an idempotent. By Theorem \ref{main1}, there exists a decomposing surface $F$ for $M$ and morphisms $A\in Mor(G,F)$, $B\in Mor(F,G)$ such that $A\circ B \cong F\times I \in Mor(F,F)$ and $B\circ A \cong M$. Since $F$ is a compact orientable surface, both $F$ and $F\times I$ embed in $\mathbb{R}^3$. Hence, $M$ embeds as $B\circ A$ in $A\circ B \circ A \circ B$. In turn, $A\circ B \circ A \circ B$ can be identified with any embedding of $F\times I$ in $\mathbb{R}^3$. So, $M$ embeds in $\mathbb{R}^3$.
\end{proof}

In light of Proposition \ref{AllEmbed}, we instead focus on which idempotents embed effectively. 

\begin{definition}
An idempotent 3-manifold $M\in Mor(F,F)$ has an \emph{effective} embedding into $\mathbb{R}^3$ if the image of $M$ in $\mathbb{R}^3$ can be surgered along the decomposing surface $F$ corresponding to $M=M\circ M$ to produce two embeddings of $M$, denoted $M_1$ and $M_2$, such that there is an embedded $2$-sphere in $\mathbb{R}^3$ separating $M_1$ from $M_2$.
\end{definition}

The following is the proof of Theorem \ref{main2}.

\begin{figure}[h!]
\begin{picture}(190,195)
\put(1,1){\includegraphics[scale=.5]{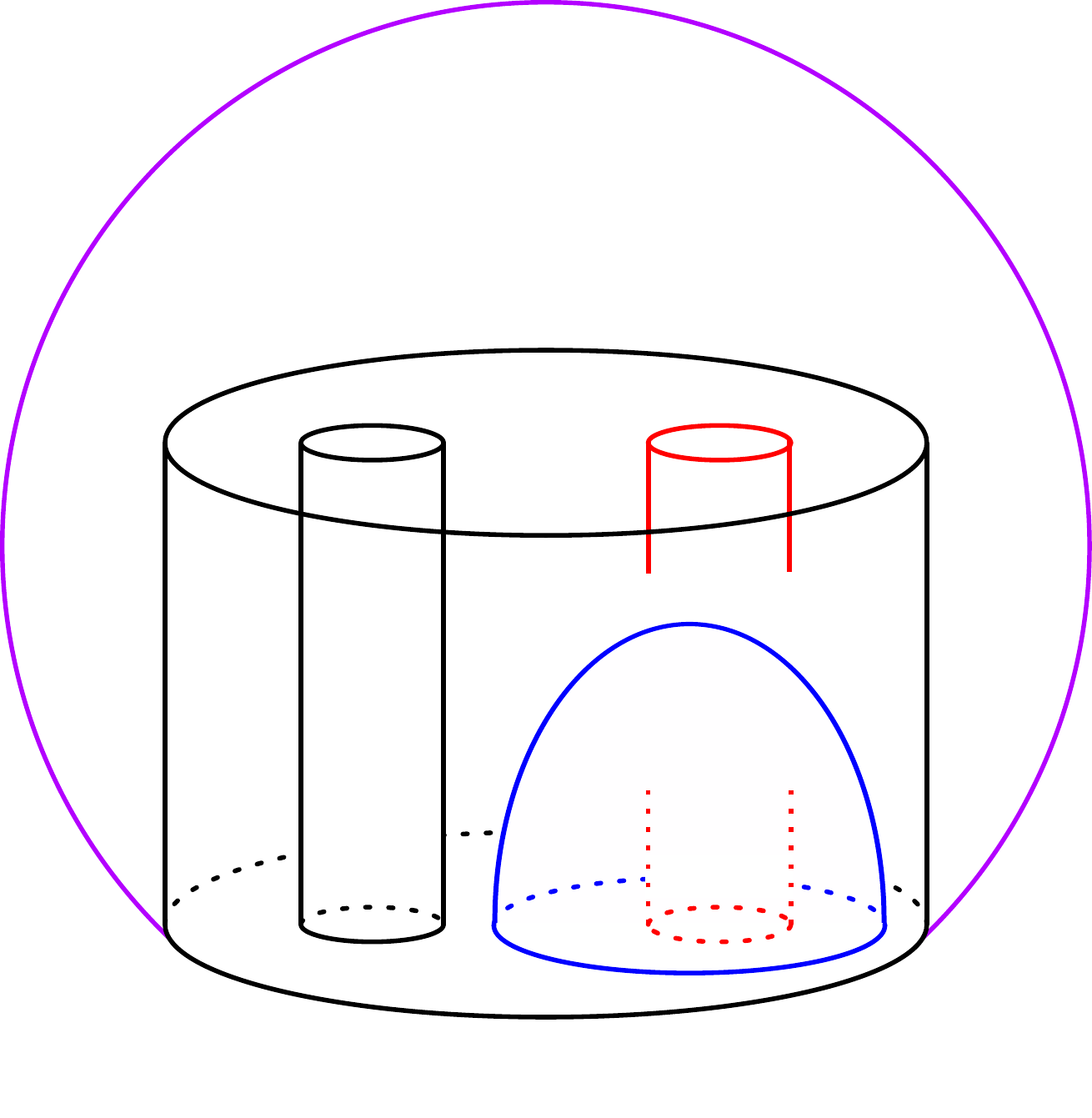}}
\put(140,113){$A$}
\put(167,167){$S$}
\put(144,74){$D$}
\put(92,155){$B_1$}
\put(115,65){$B_2$}
\end{picture}
\caption{$M_1$ is embedded in the 3-ball $B$ with the compressing disk $D$ for $\partial_+ M_1$ in blue and the annulus $A$ in red.}
\label{effective}
\end{figure}
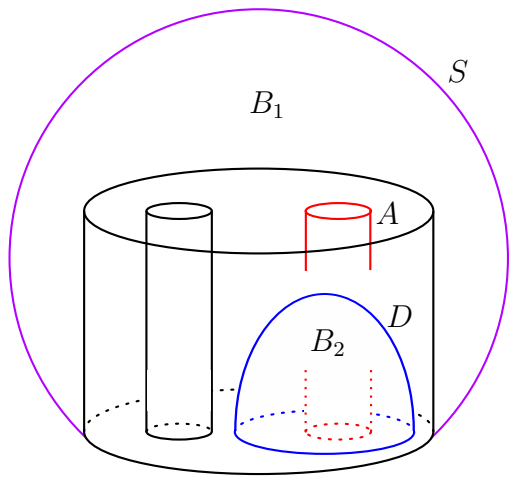

\begin{proof}
If $G$ is planar and $M\in Mor(G,G)$ is trivial, then an effective embedding of $M$ in $\mathbb{R}^3$ is given by the product of an embedding of $G$ in $\mathbb{R}^2$ with an embedding of $[0,1]$ in $\mathbb{R}$. 

Suppose a connected idempotent $M\in Mor(G,G)$ has an effective embedding $f:M\rightarrow \mathbb{R}^3$.  Let $G'$ be a decomposing surface in $M$ which realizes $M=M\circ M$.  Let $N\cong G\times I$ be a fibered neighborhood of $G'$ in $M$ such that $G'=G\times \{\frac{1}{2}\}$ and $N\cap \partial M=\partial(G) \times I$ is a collection of annuli. Then we can surger $M$ along $G'$ to construct a new 3-manifold $M\setminus [G\times (0,1)]$, which has two connected components. We denote the image of these connected components in $\mathbb{R}^3$ by $M_1$ and $M_2$. Since $f$ is an effective embedding of the idempotent $M$, each of $M_1$ and $M_2$ are homeomorphic to $M$, and there exists an embedded 2-sphere $S$ in $\mathbb{R}^3$ such that $S$ separates $M_1$ from $M_2$.

Over all spheres in $\mathbb{R}^3$ that separate $M_1$ from $M_2$, choose $S$ to minimize the number of curves of intersection with $\partial f(N)$. If $S\cap \partial f(N)$ is empty, then $f^{-1}(S)$ is a properly embedded sphere in $N\cong G\times I$ that separates $G\times \{0\}$ from $G\times \{1\}$, which is impossible unless $G\cong S^2$. Since every smoothly embedded 2-sphere in $\mathbb{R}^3$ bounds a 3-ball to one side and $\partial f(M)$ is the disjoint union of two 2-spheres, then $f(M)$ is homeomorphic to $S^2\times I$ and $M$ is a trivial morphism. Hence, we can assume $S\cap \partial f(N)$ is a non-empty collection of curves. 

Suppose $S\cap f(N)$ is compressible in $f(N)$ with compressing disk $D$. Since $S$ is a sphere, then $S\cap f(N)$ is planar and $\partial D$ separates boundary components in $S\cap f(N)$. Surgering $S$ along $D$ produces two embedded 2-spheres in $\mathbb{R}^3$, each intersecting $\partial f(N)$ in strictly fewer curves than $S$. At least one of these two 2-spheres also separates $M_1$ from $M_2$, contradicting the minimality of $S\cap f(N)$. Hence, $S\cap f(N)$ is an incompressible surface properly embedded in a 3-manifold homeomorphic to $G\times I$ which separates $G\times \{0\}$ from $G\times \{1\}$. Since all incompressible surfaces in $(surface)\times I$ are vertical or horizontal, then $S\cap f(N)$ is horizontal and, hence, every component of $S\cap f(N)$ is properly isotopic to $f(G\times \{0\})$ in $f(N)$. Thus, $G$ is a planar surface.

In the following paragraph we identify $\mathbb{R}^3$ with its one point compactification, $S^3$ so that the sphere $S$ bounds a $3$-ball to each side. Since every component of $S\cap f(N)$ is properly isotopic to $f(G\times \{0\})$ in $f(N)$, then there is an isotopy of $S$ after which $M_1$ is embedded in a 3-ball $B$ bounded by $S$ and $M_1\cap S=\partial_+ M_1\cong G$. Suppose $\partial_+ M_1$ is compressible in $M_1$ with compressing disk $D$. Since $\partial D$ is essential in $\partial_+ M_1$ and $\partial_+ M_1$ is planar, then $\partial D$ separates the boundary components of $\partial_+ M_1$. Let $E_1$ and $E_2$ be the two disks in $S$ bounded by $\partial D$. The disk $D$ cuts $B$ into two 3-balls $B_1$ and $B_2$ such that $\partial B_1=E_1 \cup D$ and $\partial B_2=E_2 \cup D$. Since $\partial_- M_1$ is connected, then we can assume $\partial_- M_1$ is embedded in $B_1$, up to relabeling. Let $\alpha$ be a boundary component of  $\partial_+ M_1$ contained in $E_2$. The curve $\alpha$ together with some boundary component $\beta$ of $\partial_- M_1$ cobound an annulus $A\subset \partial M_1$. See Figure \ref{effective}. However, $int(B_1)$ and $int(B_2)$ induce a separation on $A$, contradicting the fact that the annulus is connected. Hence, $\partial_+ M_1$ is incompressible in $M_1$. By isotoping $S$ so that $M_2$ is embedded in the other 3-ball bounded by $S$ and $M_2\cap S=\partial_- M_2\cong G$ and repeating the above argument, we can show $\partial_- M_2$ is incompressible in $M_2$. Taken together this implies that both $\partial_+ M$ and $\partial_- M$ are incompressible in $M$. By the claim in the proof of Lemma \ref{Essential}, $M$ is a trivial morphism.

\end{proof}

\section{acknowledgements}  The first author was partially supported by NSF Grant DMS--1821254.

\bibliographystyle{plain}
\bibliography{bib}

\end{document}